\documentclass[10pt,a4paper]{amsart}
\usepackage[latin1]{inputenc}
\usepackage{amsmath}
\usepackage{amsfonts}
\usepackage{amssymb}
\usepackage{graphicx}
\usepackage{amscd}
\usepackage{mathrsfs}
\usepackage[all]{xy}
\usepackage{overpic}

\newtheorem{thm}{Theorem}[section]

\newtheorem{prop}[thm]{Proposition}

\theoremstyle{definition}

\theoremstyle{remark}

\numberwithin{equation}{section}

\newcommand{\tr}{\text{tr}}                        

\newcommand{\dbar}{\bar{\partial}}        
\newcommand{\hook}{\text{{\Large$ \lrcorner\,$}}}        
\newcommand{\J}{\mathcal{J}}                          

\newcommand{\Comp}{\mathbb{C}}       
\newcommand{\Real}{\mathbb{R}}       
\newcommand{\Ham}{\mathbb{H}}        
\newcommand{\Oct}{\mathbb{O}}         
\newcommand{\R}{\mathbb{R}}      

\newcommand{\CP}{\mbox{$\mathbb{CP}$}}        

\begin{document}

    \author{Andrew Clarke}
    \title{Minimal Surfaces in $G_2$-Manifolds}    
 \address{Laboratoire de Math{\' e}matiques Jean Leray (UMR 6629) , Universit{\' e} de Nantes,  
            2, Rue de la Houssini{\` e}re, B.P. 92208, 44322 Nantes Cedex 3, France}
    \email{andrew.clarke@univ-nantes.fr}
    
\date{\today}

\maketitle


\begin{abstract}
We consider immersions of a Riemann surface into a manifold with $G_2$-holonomy and give criteria for them to be conformal and harmonic, in terms of an associated Gauss map.
\end{abstract}

\section{Introduction}

The decomposition of $\Real^7$ as the sum $\Real^3+\Real^4$ has been utilised almost since the initial advent of $G_2$-geometry. In one form this is seen in the distinguished families of minimal submanifolds that arise from the calibration condition for the degree $3$ and $4$ forms that define the structure. In another way, this is seen in the development of a gauge theory for $G_2$-manifolds that mirrors the theories that exist in $3$ and $4$ dimensions. 

In the first case, given the differential $3$-form that determines the metric, one can consider those $3$-dimensional submanifolds $L$ that satisfy the calibration condition $\varphi|_L=\text{vol}_L$. Any such submanifold is homologically volume minimising. This is to say that, in addition to being a minimal submanifold, the volume of $L$ is minimal among all submanifolds in the same (relative) homology class. In dimensions other than $3$ and $4$ however, we can not consider submanifolds that relate to the ambient spacein this way. This is the subject of this article. We consider maps from a $2$-dimensional Riemannian manifold into a space $X$. 

The $4$-dimensional version of the principal results given here were obtained by Eells and Salamon \cite{ES}. Over an oriented Riemannian $4$-manifold $M$ one can consider the positive and negative twistor spaces $Z_+$ and $Z_-$. that fibre over $M$. These spaces can be defined in many ways but the one that is tractable in this case is to consider them to be the sphere bundles of subspaces of $\Lambda^2$ that correspond to irreducible representations of $SO(4)$. Among many other results, they prove that for any immersion of a surface $f:\Sigma\to M$, the adapted Gauss lift $\tilde{f}_-:\Sigma\to Z_-$ is $J_1$-holomorphic if and only if $f$ is conformal and harmonic. Here $J_1$ is a naturally defined almost-complex structure on $Z_-$. We prove an analogue of this result. 

In the current paper we consider a sphere fibration of a $G_2$-manifold $\pi:Z_7\to X$ where in this case $Z_7$ is the unit tangent sphere bundle. For any immersion of an oriented surface $F:\Sigma\to X$, the vector cross product on each tangent space defines a lifting $\tilde{f}_7:\Sigma\to Z_7$. We show that the space $Z_7$ admits a distribution $\mathcal{E}$ with Hermitian structure defined in terms of the $G_2$-geometry of the base. For any immersion $f$ into the base, the lifting $\tilde{f}_7$ of $f$ is everywhere tangent to the distribution $\mathcal{E}$.  If $\Sigma$ has a given conformal structure, one can give geometric criteria for the derivative $\tilde{f}_{7\, *}:T_\Sigma\to\mathcal{E}$ to be complex linear. The principal hypothesis is that $f$ be conformal and harmonic, the geometric meaning of which is clear. The second hypothesis is that the second fundamental form takes values orthogonal to the unit vector $\eta$ that defines the lifting $\tilde{f}_7$. In Section \ref{sec:adapted} this condition is interpreted in terms of the holomorphicity of a vector bundle associated to the immersion. In the final section, the results are interpreted on a number of important $G_2$-manifolds.

This analogue can be understood by considering a recent theorem of Verbitsky \cite{V}. One of the celebrated theorems of twistor theory states that $Z_-$ admits an almost-complex structure $J_2$ and this structure is integrable if and only if a component of the Riemannian curvature tensor of $M$ vanishes. Verbitsky proves that if $(X,\varphi)$ is a $7$-manifold with $G_2$-structure, the manifold $Z_7$ admits an almost $CR$-structure that is integrable if and only if the $G_2$-structure is torsion free. This can be interpreted in terms of some components of the Riemannian curvature tensor vanishing. 

The use of Gauss maps to study minimal submanifolds is a well established technique. In addition to the work of Eells and Salamon, it was used by Bryant and Gauduchon \cite{B,G} to study super-minimal immersions. This is another notion that could be also studied in the $G_2$ context, although we leave this to another occasion. Gauss maps can be used to study other curvature conditions on submanifolds. We take the opportunity to mention the work of Labourie and Smith \cite{L,S} as examples.

This project was largely undertaken at {\'E}cole Polytechnique and the Instituto Nacional de Matem{\'a}tica Pura e Applicada. The author would like to thank these institutions for their hospitality and excellent working environment.

\section{Linear Algebra of $2$-planes in $\R^7$}\label{sec:planes}

$G_2$-geometry is defined in terms of a three-form on $\R^7$. The three form is defined so as to encode the multiplicative structure on $\Oct$ and the vector cross product structure on $\R^7=\text{Im}\Oct$. Following Karigiannis \cite{K1} and Harvey and Lawson \cite{HL} we define
\begin{eqnarray*}
\varphi_0(x,y,z)=\langle x,yz\rangle =\langle x,y\times z\rangle.
\end{eqnarray*}
The cross product $\times$ is skew symemtric and satisfies
\begin{eqnarray*}
|y\times z| = |y\wedge z|
\end{eqnarray*}
where the norm on the right is of the simple 2-vector in $\Lambda_2\R^7$. With respect to the standard orthonormal basis $\{e_1,\ldots,e_7\}$ for $\R^7$ can be expressed as
\begin{eqnarray*}
\varphi_0&=&\varepsilon^{123}+\varepsilon^i\wedge \eta^-_i\\
&=& \varepsilon^{123}+\varepsilon^1\wedge(\varepsilon^{45}-\varepsilon^{67})+ \varepsilon^2\wedge(\varepsilon^{46}-\varepsilon^{75}) +\varepsilon^3\wedge(\varepsilon^{47}-\varepsilon^{56})
\end{eqnarray*}
where $\{ \varepsilon^i\}$ form the dual basis of covectors and $\varepsilon^{ij}=\varepsilon^i\wedge\varepsilon^j$.
We note that  $\{\eta^-_i\}$ form a basis of anti-self-dual $2$-forms on $\R^4$ (see \cite{K1}). 

Given an oriented 2-plane $\xi=\text{span}(y,z)$ of $\R^7$ where $y,z$ form an oriented orthonormal basis for $\xi$, we can define a vector $\eta_\xi=y\times z\in S^6$. Due to the skew symmetry of the cross product, this is independent of the oriented orthonormal basis for $\xi$. $\eta(y\wedge z)$ is uniquely determined to be that unit vector such that $y\wedge z\wedge \eta$ is an associative $3$-plane.

 For example, for $\xi=e_1\wedge e_2$, we have
\begin{eqnarray*}
\eta_\xi = ((e_1\wedge e_2)\hook\varphi_0)^\#.
\end{eqnarray*}

In fact, this example effectively describes the definition and much of the pointwise geometry for an arbitrary $2$-plane. In \cite{HL} it is shown that $G_2$ acts transitively on the Stiefel manifold $V_{2,7}$ of ordered orthonormal pairs of vectors in $\R^7$. That is, up to the action of $G_2$, we can say that an arbitrary $2$-plane $\xi$ is spanned by $\{e_1,e_2\}$ where $\varphi_0$ has the form 
\begin{eqnarray}
\varphi_0=\varepsilon^{123}+\varepsilon^i\wedge \eta^-_i.\label{phi34}
\end{eqnarray}

Given an oriented $2$-plane $\xi\subseteq \Real^7$ we obtain a unit vector $\eta$ orthogonal to $\xi$. This vector determines the complex structure arising from the orientation and inner product by the definition
\begin{eqnarray*}
\langle \mathcal{J}v,w\rangle &=& \varphi_0(\eta,v,w)\\
&=& \langle \eta\times v,w\rangle.
\end{eqnarray*}
That $\mathcal{J}$ satisfies $\mathcal{J}^2=-1$ follows from contractions and restrictions of $\varphi$ to $\xi$. The corresponding  fundamental $2$-form is given by $\omega=(\eta\hook\varphi)|_{\xi}$. Similarly, left cross multiplication by $\eta$ preserves the subspace $W=\text{span}\{e_4,e_5,e_6,e_7\}$. The fundamental form on $W$ is given by 
\begin{eqnarray*}
\omega=(\eta\hook\varphi)|_W.
\end{eqnarray*}



In particular, this shows explicitly how $S^6\subseteq\R^7$ is endowed with a $G_2$-invariant almost complex structure. The tangent space to $S^6$ at $\eta$ is the space $E_\eta=\langle\eta\rangle^\perp$. As we have seen, left multiplication by $\eta$ on this vector space defines the structure. In the expression $S^6=G_2/SU(3)$ of $S^6$ as a homogeneous space, this is the almost complex structure given by the $SU(3)$-structure.

The action of $SO(7)$ on $\Lambda^2\R^7$ is irreducible but when restricted to $G_2$ this is no longer true. There is a decomposition
\begin{eqnarray*}
\Lambda^2\R^7 = \Lambda^2_7\oplus \Lambda^2_{14}.
\end{eqnarray*}
According to the identification $\Lambda^2\R^7\cong \mathfrak{so}(7)$ by raising an index the component $\Lambda^2_{14}$ is identified with $\mathfrak{g}_2$. We also have identifications 
\begin{eqnarray}
\Lambda^2_{14} &=& \ker\{ (*\varphi_0)\wedge:\Lambda^2\to \Lambda^6\}\nonumber\\
\Lambda^2_{7} &=& \{ \frac{1}{3} X\hook\varphi_0\ ;\ X\in \R^7\} \label{eqn:L7isomR7}.
\end{eqnarray}
More important for our interests at hand are that $\Lambda^2_7$ and $\Lambda^2_{14}$ are eigenspaces of the operation on $\Lambda^2\R^7$
\begin{eqnarray*}
\alpha\mapsto *(\varphi_0\wedge\alpha).
\end{eqnarray*}
Using the convention that $\varphi_0=\varepsilon^{123}+\varepsilon^i\wedge\eta^-_i$ we have
\begin{eqnarray*}
\Lambda^2_7 &=& \{\alpha\ ;\ *(\varphi_0\wedge\alpha)=-2\alpha\}\\
\Lambda^2_{14}&=&\{ \alpha\ ;\ *(\varphi_0\wedge\alpha)=\alpha\}.
\end{eqnarray*}
The projections from $\Lambda^2$ to the two factors can be given 
\begin{eqnarray*}
\pi_7 &=&\frac{1}{3}(Id-*\circ\varphi_0\wedge)\\
\pi_{14}&=& \frac{1}{3}(2Id +*\circ \varphi_0\wedge).
\end{eqnarray*}
These calculations are expressed in terms of exterior covectors, which are typically the natural objects to study. In the case at hand though, for the embeddings of Grassmannians and for the study of Gauss maps, we consider the algebra of exterior vectors, denoted $\Lambda_p\Real^n$. These spaces have identical decomposition to those above and only differ notationally in the presence of a metric. Notationally though, they are simpler in this case.

The Grassmannian $G(2,7)$ of oriented $2$-planes in $\R^7$ naturally embeds in $\Lambda_2\R^7$. A plane $\xi$ with oriented orthonormal basis $\{e_1,e_2\}$ can be identified with $e_1\wedge e_2\in \Lambda_2\R^7$. As noted above, $G_2$ acts transitively on $G(2,7)$, with isotropy group $U(2)$, and this is compatible with the representation of $G_2$ on $\Lambda_2\R^7$ induced from that on $\R^7$. 

The action of these projections on the Grassmannian in this space can also be understood. The image $\pi_7(G(2,7))$ is an orbit of $G_2$ on $\Lambda_{2,7}$ and can be identified as $S^6\subseteq \R^7$. In fact,
\begin{eqnarray}
\pi_7(e_1\wedge e_2)=\frac{1}{3}(e_{12}+e_{47}-e_{56})=\frac{1}{3}(e_3\hook\varphi_0)^\# \label{eqn:pi7is1}
\end{eqnarray}
where $e_3=\eta=e_1\times e_2$ is the vector orthogonal to $\xi=e_1\wedge e_2$ that is considered above. In terms of homogeneous spaces,
\begin{eqnarray*}
G(2,7)=G_2/U(2) \to G_2/SU(3) = S^6
\end{eqnarray*}
with fibre $SU(3)/U(2) = \CP^2$ with the fibre over $\eta\in S^6$ being the set of $\J_\eta$-complex lines in $\langle\eta\rangle^\perp$ (with orientation). On the other hand, $\pi_{14}:\Lambda_2\R^7\to \Lambda_{2,14}$ is injective when restricted to $G(2,7)$ and defines an equivariant embedding of $G(2,7)$ into $\mathfrak{g}_2$. In particular, this realises $G(2,7)$ as an adjoint orbit of $G_2$.

According to the standard embedding $G(p,n)\subseteq \Lambda_p\R^n$, the tangent space to $G(p,n)$ at $\xi$ is generated by exterior elements of the form
\begin{eqnarray*}
a\wedge(b\hook\xi)  \ \ \ \ \ \text{for}\ b\in\xi\ \text{and}\ a\perp\xi.
\end{eqnarray*}
For $G(2,7)\subseteq \Lambda_2\R^7$ and for $\xi=e_1\wedge e_2$ we can then see that
\begin{eqnarray*}
T_\xi G(2,7)=\{ V_1\wedge e_1 +V_2\wedge e_2 +e_3\wedge (ae_1+be_2)\}
\end{eqnarray*}
for $V_1,V_2\in W=\langle e_1,e_2,e_3\rangle^\perp$ and $a,b\in\R^7$. A $G_2$-invariant almost complex structure can be defined on $G(2,7)$ by denoting 
\begin{eqnarray*}
T_\xi^{1,0}=W^{1,0}\wedge e_1\oplus W^{1,0}\wedge e_2\oplus e_3\wedge\xi^{1,0}\subseteq T_\xi\otimes \Comp
\end{eqnarray*}
where the subspaces $\xi$ and $W$ have the almost complex structures given, as above, by left multiplication by $e_3=\eta(\xi)$. When we express the Grassmannian as a homogeneous space $G(2,7)/U(2)$, it can be seen that the structure given above is the one compatible with the isotropy action of $U(2)$.

We need to consider the projections
\begin{eqnarray*} 
\pi_7 &:& G(2,7) \to S^6\\
\pi_7 &:& T_\xi G(2,7)\to T_\eta S^6\ \ \ \eta=\pi_7(\xi).
\end{eqnarray*}
It can easily be seen that $\pi_7$ is not complex linear, but in fact $\pi_7(e_3\wedge \xi^{1,0})\subseteq E^{1,0}$ and $\pi_7(W^{1,0}\wedge e_i)\subseteq E^{0,1}$. This is to say that $\pi_7$ is linear, and anti-linear on the respective subspaces. These statements can be made more precise. We can furthermore see that 
\begin{eqnarray*}
\ker\pi_7 =W^{1,0}\otimes \xi^{0,1} + W^{0,1}\otimes \xi^{1,0}
\end{eqnarray*}
and that the following maps are bijections
\begin{eqnarray*}
\pi_7 &:& W^{1,0}\otimes \xi^{1,0} + e_3\otimes \xi^{0,1} \to E^{0,1}\\
\pi_7 &:& W^{0,1}\otimes \xi^{0,1} + e_3 \otimes \xi^{1,0} \to E^{1,0}.
\end{eqnarray*}
In particular, for $\alpha\in T_\xi G(2,7)\otimes \Comp\subseteq \Lambda_2\Real^7\otimes \Comp$, and for $\varepsilon=e_1-ie_2\in \xi^{1,0}$, the projection $\pi_7(\alpha)\in E^{1,0}$ if and only if the component of $\varepsilon\hook\alpha$ in $W^{1,0}$ is equal to zero and the component of $\bar{\varepsilon}\hook\alpha$ in $\langle e_3\rangle$ is equal to zero.

\section{Manifolds with $G_2$ Holonomy}

We continue the above discussion and consider these constructions to be the pointwise description of our geometry.

Let $(X,\varphi,g)$ be a Riemannian $7$-manifold with holonomy group contained in $G_2$. This is equivalent to the following. $X$ possesses a smooth $3$-form $\varphi$ such that at each point there exists a basis $\{e_i\}$ of $T_X$ with dual coframe of $1$-forms $\{\varepsilon^i\}$ such that $\varphi$ is expressed at $p$ as 
\begin{eqnarray*}
\varphi=\varepsilon^{123}+\varepsilon^i\wedge\eta^-_i,
\end{eqnarray*}
ie., $\varphi$ is pointwise equivalent to $\varphi_0$ on $\R^7$. Just as on $\R^7$, $\varphi$ determines an orientation $d\mu$ and metric $g_\varphi$ on $X$. The holonomy group of $(X,g_\varphi)$ is contained in $G_2$ if and only if $\varphi$ is parallel with respect to the metric that it determines. 

A consequence of this definition is that on each tangent space there exists a vector cross product isomorphic to the one on $\R^7$. The exterior bundles decompose into components isomorphic to those on $\R^7$. The Grassmann bundle can at each point be projected onto the $7$ and $14$-dimensional components of $\Lambda^2T^*_X$. These properties hold for any $7$-dimensional  manifold for which the structure group reduces to $G_2$. If $Hol(g_\varphi)\subset G_2$ though, the Levi-Civita connection preserves the decomposition and acts essentially identically on isomorphic components. 

We consider the Grassmann bundle $G_2(X)$ of oriented $2$-planes tangent to $X$. This space fibres over $X$ with each fibre diffeomorphic to $G(2,7)$. By the same construction given earlier, the bundle $G_2(X) $ is naturally contained in the vector bundle $\Lambda_2T_X$ of exterior $2$-vectors. In the same fashion as for differential forms, this space splits into $7$ and $14$-dimensional components; $\Lambda_{2,7}$, $\Lambda_{2,14}$. The component $\Lambda_{2,7}$ is naturally isomorphic to the tangent bundle $T_X$. 

Again following the previous section we can consider the projection $\pi_7$ from $\Lambda_2$ to $\Lambda_{2,7}$. As previously noted, the restriction of $\pi_7$ to $G(2,7)\subseteq\Lambda_2\Real^7$ defines a fibration of $G(2,7)$ over the six-sphere. We can consider this to be the model at each point and conclude that the restriction of $\pi_7$ to $G_2(X)$ defines a map $\pi_7:G_2(X)\to Z_7$ onto the unit sphere bundle of the tangent bundle to $X$. 

We now recall some definitions of connections on geometric fibre bundles. Let $X$ be an $n$-dimensional Riemannian manifold and let $G_p(X)$ be the Grassmann bundle of oriented $p$-planes. The metric on $X$ defines a splitting of the tangent bundle of $G_p(X)$. At $\xi\in G_p(X)$, $T_\xi G_p(X)=V_\xi + H_\xi$ where the vertical space $V_\xi=\ker(\pi_*)$ and the horizontal space $H_\xi$ is generated by infinitesimal parallel translates of $\xi$ along paths in $X$ through $x=\pi(\xi)$. $\pi_*$ defines an isomorphism $\pi_*:H_\xi\to T_xX$. For $v\in T_xX$ the horizontal lift of $v$ to $\xi$ is given as follows. Along a path in $X$ tangent to $v$ we extend $\xi$ to a family of $p$-planes and define the horizontal lift of $v$ at $\xi$ to be 
\begin{eqnarray}\label{eqn:horizlift}
v_\xi^H={\dot{\xi}}-(\nabla_v\xi)^V
\end{eqnarray}
where the final term is the covariant derivative of $\xi$ along the path in $X$, with respect to the Levi-Civita connection. This definition is independent of the extension of $\xi$ that is chosen. It is depends only on the point $\xi\in G_p(X)$ and the vector $v\in T_xX$. 

We apply this in particular to the case $p=1$ for the unit sphere bundle $Z_7=U_1(T_X)$ of the tangent bundle to $X$. In this case, at a point $\eta\in Z_7$ the vertical space $V_\eta$ can be identified with the subspace $E_\eta=\langle \eta\rangle^\perp\subseteq T_{x,X}$. The manifold $Z_7$ canonically admits a contact structure. Given the projection map $\pi:Z_7\to X$, we define 
\begin{eqnarray*}
\mathcal{E}_\eta = \pi_*^{-1}(\langle\eta\rangle^\perp)\subseteq T_\eta Z_7.
\end{eqnarray*} 
The tangent space to $Z_7$ splits as a sum of vertical and horizontal components $T_\eta Z_7=V_\eta + H_\eta$. The projection $\pi$ is an isomorphism to $T_X$ when restricted to the horizontal subspace. We therefore have a decomposition
of $\mathcal{E}$, 
\begin{eqnarray*}
\mathcal{E}_\eta= E^V_\eta + E^H_\eta
\end{eqnarray*}
where both $E^V_\eta=V_\eta $ and $E^H_\eta =\pi_*^{-1}(E_\eta)\subseteq H_\eta$ are canonically isomorphic to the subspace $E_\eta=\langle \eta\rangle^\perp\subseteq T_{x,X}$. As we have seen in the previous section, for any unit vector $\eta\in T_{x,X}$ the subspace $E_\eta$ can be given an almost complex structure by defining
\begin{eqnarray*}
\mathcal{J}:E_\eta&\to &E_\eta\\
\langle \mathcal{J}v,w\rangle &=&\varphi(\eta,v,w).
\end{eqnarray*}
Using the above identifications we can consider $E^V$ and $E^H$ as having almost complex structures. We can in turn define almost complex structures on the distribution $\mathcal{E}$ on $Z_7$. The structure in which we will be interested for the moment can be defined as follows. We take
\begin{eqnarray*}
\mathcal{E}^{1,0} = (E^V)^{1,0}+ (E^H)^{1,0}\subseteq \mathcal{E}\otimes \Comp.
\end{eqnarray*}
This is to say that we take $\mathcal{J}$ on both the horizontal and vertical summands of $\mathcal{E}$.

It should be noted that this structure differs from the almost complex structure usually given on the canonical contact distribution on the unit tangent bundle to a Riemannian manifold. In that case, the automorphism interchanges the horizontal and vertical components. In many cases, geometric conditions can often be given for the distribution with complex struture to define an integrable CR structure. This does not apply however in the case at hand however. The almost complex structure on $\mathcal{E}$ restricts to the vertical subspace to define the standard non-integrable $G_2$-invariant almost complex structure on the $S^6$ fibres.

\section{Harmonic Maps and Minimal Surfaces}

Let $(M,g)$ and $(N,h)$ be Riemannian manifolds and $f:M\to N$ a smooth map between them. The derivative between them can be considered a section 
\begin{eqnarray*}
df\in\Gamma(M,T^*_M\otimes f^*T_N)
\end{eqnarray*}
of a Riemannian vector bundle over $M$. The \emph{tension field} $\tau_f$ of $f$ is defined to be the divergence of this section
\begin{eqnarray*}
\tau_f=\tr(\nabla df)\in \Gamma(M,f^*T_M).
\end{eqnarray*}
$f$ is said to be a \emph{harmonic} map if $\tau_f=0$. If $f$ is an isometric immersion, the Levi-Civita connection of $M$ is obtained as the tengential projection of the Levi-Civita connection on $N$. The torsion of $f$ can then be related to the second fundamental form of the immersed submanifold. Let $M\subseteq N$ be a smooth submanifold. The second fundamental form $B$ of $M$ is a bilinear form on $T_M$ with values in the normal bundle to $M$. 
\begin{eqnarray*}
B(v,w)&=&(\nabla_vw)^N\\
B&\in &\Gamma(S^2(T^*_M)\otimes N_M).
\end{eqnarray*}
The mean curvature $H$ is the trace of $B$. If $H$ vanishes, $M$ is said to be a minimal submanifold. It can easily be seen that if $f:M\to N$ is an isometric immersion, $f$ is a harmonic map if and only if the image $f(M)$ is a minimal submanifold. 
This can be strengthened in $M$ is $2$-dimensional. The condition that $f$ be harmonic is, for a fixed metric on $N$, dependent only on the conformal structure on $M$. We can then conclude that the image of a conformal immersion is minimal if and only if the immersion is harmonic. 

For a smooth manifold $X$ of dimension $n$ we can consider the fibration $\pi:G_p(X)\to X$ given by the set of all oriented tangent $p$-planes to $X$. At each point $x\in X$ the fibre over $x$ consists of the $p$-dimensional planes contained in $T_xX$ with a fixed orientation. This fibre is diffeomorphic to the Grassmannian $G(p,n)$ of oriented $p$-planes in $\Real^n$. 

Let $M$ be a $p$-dimensional oriented manifold and $f:M^p\to X^n$ a smooth immersion. The Gauss lift of $f$ is defined to be 
\begin{eqnarray*}
{\tilde f}:M&\to & G_p(X)\\
x& \mapsto &f_*(T_xM)
\end{eqnarray*}
 If $X$ is a Riemannian manifold we can use the embedding $G(p,n)\to \Lambda_p\Real^n$ to describe the Gauss lift using exterior algebra. let $f: M\to X$ be an isometric immersion and at $x\in M$ let $\{e_1,\ldots,e_p\}$ be an oriented orthonormal basis for $T_xM$. Then ${\tilde f}$ is given by 
\begin{eqnarray*}
{\tilde f}(x)=e_1\wedge \cdots \wedge e_p.
\end{eqnarray*}
In the case where $M=\Sigma$ is of real dimension $2$ and if $f:\Sigma\to X$ is an isometric immersion, we can consider the orthonormal basis $\{e_1,e_2\}$ such that $\varepsilon=e_1-ie_2$ is of type $(1,0)$ with respect to the induced complex structure on $\Sigma$. Then, $e_1\wedge  e_2=-i/2\varepsilon\wedge {\bar \varepsilon}$.

We can extend the observation given in Equation \ref{eqn:horizlift} to consider the derivative of the Gauss lift. Let $f:M\to X$ be an immersion and $\tilde{f}:M\to G_p(X)$ be a lifting of $f$ that satisfies $\pi\circ\tilde{f}=f$. Then $\tilde{f}$ defines a section $\sigma$ of the bundle $f^*(\Lambda_pT_X)\to M$. Then for $v\in T_M$, the horizontal and vertical parts of $\tilde{f}_*(v)$ are given by 
\begin{eqnarray*}
\tilde{f}_*(v)=(\tilde{f}_*(v))^H +((f^*\nabla)_v\sigma)^V
\end{eqnarray*}
where $f^*\nabla$ is the pull-back of the Levi-Civita connection on $X$.

\section{Complex Linear Gauss Lifts}

We now suppose that the manifold $X$ is $7$-dimensional and equipped with a metric of holonomy $G_2$, given by a $3$-form $\varphi$. The integrability of the $G_2$-structure is required in this case because we ask that the connection preserve the decomposition $\Lambda_2=\Lambda_{2,7}+\Lambda_{2,14}$, and that the connection on the $7$-dimensional component can be identified with the connection acting on vector fields. 

Let $f:\Sigma\to X$ be an immersion. 
Then, in addition to the usual Gauss map $\tilde{f}:\Sigma\to G_2(X)$, we can consider the projection to the $7$-dimensional component $\tilde{f}_7=\pi_7\circ\tilde{f}:\Sigma\to Z_7$. Then an elementary initial observation is the following.
\begin{prop}
For any immersion $f:\Sigma\to X$, $\tilde{f}_7(\Sigma)$ is an integral submanifold of the distribution $\mathcal{E}$ on $Z_7$. 
\end{prop}
\proof{ 
For $f:\Sigma\to X$, $\pi(\tilde{f}_7(p))=f(p)$ so $\pi_*(\tilde{f}_{7*}(v))=f_*(v)$ which is perpendicular to the vector $\eta$ that defines both the lift $\tilde{f}_7$ and the distribution $\mathcal{E}$. It clearly follows that $\tilde{f}_{7*}(v)$ is tangent to $\mathcal{E}$. 
}

That is, for any immersion $f:\Sigma\to X$, the linearisation of the $7$-dimensional Gauss lift defines a linear map $\tilde{f}_{7*}:T_\Sigma\to\mathcal{E}$. If $\Sigma$ is endowed with a metric $h$, we look for conditions for when this map is complex linear. 

The horizontal and vertical components of $\mathcal{E}$ are each complex vector spaces with respect to the structure. It is sufficient that  these components of the Gauss lift be complex linear. To begin with, we consider the horizontal component. 

\begin{prop} Let $(\Sigma,h)$ be a Riemannian surface and $f:\Sigma\to X$ an immersion. 
The horizontal component of the $7$-dimensional Gauss lift is complex linear if and only if $f$ is a conformal map.
\end{prop}
\proof{ This is an immediate consequence of the definition of the the Hermitian structure on the horizontal component of Gauss lift. In particular, $\pi_*$ defines an isomorphism between $E_\eta^H$ and $E_\eta=\langle\eta\rangle^\perp=\xi\oplus W$ for $\xi=f_*(T_\Sigma)\subseteq T_X$. Vector multiplication by $\eta$ defines the Hermitian structures on $\xi$ and $W$. The map $f$ is conformal if and only if the structure induced on $\xi$ is the same as that coming from the metric $h$ on $\Sigma$. 
}


As noted above, the derivative of the $7$-dimensional Gauss lift of the map $f$ has a horizontal and a vertical component. The map is complex linear if each is complex linear separately. In the above proposition it is seen that the horizontal component of $\tilde{f}_{7*}$  sends $T^{1,0}_\Sigma$ to $(E^H)^{1,0}$ if and only if $f$ is conformal. The vertical component is dealt with in the following result.

\begin{thm} \label{thm:holom-min}
Let $f:\Sigma\to X$ be an immersion such that the second fundamental form of $f(\Sigma)$ takes values orthogonal to the vector $\eta=\eta(\xi)$. Then $\tilde{f}_7:\Sigma\to Z_7$ is holomorphic with respect to the Hermitian structure on the distribution $\mathcal{E}$ if and only if $f:\Sigma\to X$ is conformal and harmonic. 
\end{thm}
\proof{From the previous proposition we see that $f$ must be a conformal map. We will suppose that the map is in fact an isometric immersion. The conclusions that we will obtain will be conformally invariant though, so this assumption is allowed. The vertical part of the derivative is given by
\begin{eqnarray*}
(\tilde{f}_{7*}(v))^V=\nabla_v\eta
\end{eqnarray*}
where $\nabla$ is the Levi-Civita connection on $X$. Given also the choice for the Hermitian structure on $\mathcal{E}$, $\tilde{f}_{7*}$ is complex linear if and only if 
\begin{eqnarray*}
\nabla_\varepsilon\eta=\nabla_\varepsilon\pi_7(\xi)=\pi_7(\nabla_\varepsilon\xi)
\end{eqnarray*}
takes values in the vector space $E^{1,0}_\eta\subseteq\mathcal{E}^{1,0}$ whenever $\varepsilon\in T^{1,0}_\Sigma$. We have however that
\begin{eqnarray*}
\nabla_\varepsilon\xi = \frac{-i}{2}\nabla_\varepsilon(\varepsilon\wedge\bar{\varepsilon})=\frac{-i}{2}(\nabla_\varepsilon\varepsilon) \wedge \bar{\varepsilon} + \frac{-i}{2}\varepsilon\wedge\nabla_\varepsilon\bar{\varepsilon}. 
\end{eqnarray*}
If the immersion $f:\Sigma\to X$ is an isometric immersion the Levi-Civita connection on $\Sigma$ is given by the tangential component of $\nabla$. At a point $p\in\Sigma$ we can take a hermitian frame $\{\varepsilon\}$ near $p$ such that $(\nabla\varepsilon)^T(p)=(\nabla\bar{\varepsilon})^T(p)=0$. The expressions $\nabla_\varepsilon\varepsilon$ and $\nabla_\varepsilon\bar{\varepsilon}$ are therefore (complex) normal-valued at $p$. 

We recall the observation made at the end of Section \ref{sec:planes}. A vector $\alpha\in T_\xi G(2,7)\otimes\Comp$ is mapped by $\pi_7$ to $E^{1,0}$ if and only if the component of $\varepsilon\hook\alpha$ in $W^{1,0}$ is equal to zero and the component of $\bar{\varepsilon}\hook\alpha$ in the span of $e_3$ is zero. We apply this to $\alpha=\nabla_\varepsilon(\varepsilon\wedge\bar{\varepsilon})$. Suppose that  the second fundamental form of the immersed submanifold takes values orthogonal to the vector $\eta=e_3$. Then,
\begin{eqnarray*}
\bar{\varepsilon}\hook\nabla_\varepsilon(\varepsilon\wedge\bar{\varepsilon}) &=& \nabla_\varepsilon \varepsilon\\
&=& B(e_1,e_1)-B(e_2,e_2) -2i B(e_1,e_2).
\end{eqnarray*}
which has no component in the space of $e_3$.  For the same reason, if the second fundamental form $B$ only takes values orthogonal to $e_3$, 
\begin{eqnarray*}
\varepsilon\hook\nabla_\varepsilon(\varepsilon\wedge\bar{\varepsilon}) &=&\nabla_\varepsilon\bar{\varepsilon}\\
&=& B(e_1,e_1)+B(e_2,e_2)=H
\end{eqnarray*}
 is orthogonal to $e_3$ and so lies in $W^{1,0}+W^{0,1}$. If the  immersed submanifold $f(\Sigma)$ is minimal, $H=0$ then clearly the component in $W^{1,0}$ vanishes and $\tilde{f}_{7*}$ is complex linear. 

Conversely, if the component in $W^{1,0}$ vanishes, then $H\in W^{0,1}$. It is a purely real vector  and so must vanish if it is a $-i$ eigenvector of $\mathcal{J}$. 

To conclude the result we recall that the image of an immersion of a Riemann surface is minimal if and only if the map is conformal and harmonic. Lastly, we observe that the condition of the harmonicity of the map is invariant under conformal change. This shows that the hypothesis of $f$ being an isometric immersion was not a restriction.
}

\section{Adapted Surfaces in a $G_2$-manifold.}\label{sec:adapted}

Suppose that $\Sigma\subseteq X$ is a $2$-dimensional submanifold. Then the unit vector field $\eta$ along $\Sigma$ allows us to consider the bundle $W$ over $\Sigma$  defined as the sub-bundle of the normal bundle that is orthogonal to $\eta$. This is a Hermitian vector bundle of rank $2$. The complex structure on $W$ is given by cross product with $\eta$. That is,
\begin{eqnarray*}
\langle \J v,w\rangle = \langle \eta\times v, w\rangle =\varphi(\eta,v,w).
\end{eqnarray*}

\begin{prop} 
The connection $\nabla^W$ induced from the Levi-Civita connection on $X$, is Hermitian with respect to the complex structure $\J$ if and only if the second fundamental form $B$ takes values in $W$. 
\end{prop}
\proof{
The differential form $\varphi$ being parallel gives that
\begin{eqnarray*}
(\nabla^X\varphi)(\eta,v,w)=0= d(\varphi(\eta,v,w))-\varphi(\nabla^X\eta,v,w)-\varphi(\eta,\nabla^Xv,w)-\varphi(\eta,v,\nabla^Xw).
\end{eqnarray*}
Here we note that $\varphi(\eta,\nabla^Xv,w)=\varphi(\eta,\nabla^Wv,w)$ because the $3$-dimensional subspace $T_\Sigma+\eta$ is closed under the cross product. Hence,
\begin{eqnarray*}
\varphi(\eta,(\nabla v)^T,w)=0
\end{eqnarray*}
since $w\in W$. For a similar reason, that being that $W$ is at each point a coassociative subspace, we have that
\begin{eqnarray*}
\varphi(\nabla^X\eta,v,w)=\varphi((\nabla\eta)^T,v,w)=-\langle A^\eta\times v,w\rangle
\end{eqnarray*}
where $A^\eta$ is the shape operator (or second fundamental form) for the normal vector $\eta$, considered as an endomorphism of $T_\Sigma$.  Therefore,
\begin{eqnarray*}
\langle (\nabla^W\J)v,w\rangle &=& \langle \nabla^W(\J v),w\rangle -\langle \J(\nabla^W v),w\rangle\\
&=& d(\langle \J v,w\rangle)-\langle \J v,\nabla^W w\rangle - \langle \J\nabla^W v,w\rangle\\
&=& -\langle A^\eta\times v,w\rangle.
\end{eqnarray*}
Then, the final line in this equation vanishes for all $v,w\in W$ if and only if $A^\eta=0$. This is equivalent to the condition that $B$ takes values orthogonal to $\eta$. 
}

With this in mind, the connection $\nabla^W$ defined on $W$ is $\J$-linear and so extends to $W^{1,0}$. The $(0,1)$-part of the connection satisfies $(\nabla^{0,1})^2=0$ because the base is only of one complex dimension. This implies that $\dbar=\nabla^{0,1}$ defines an analytic structure on the bundle $W$.  Sections in the kernel of this operator are the holomorphic sections.






The next statement that we make is a relationship between the vector bundle $W$ and the intrinsic geometry of $\Sigma$. Specifically we consider the map 
\begin{eqnarray*}
F:T_X&\to& \Lambda^2T^*_X\\
Y&\to& Y\hook\varphi.
\end{eqnarray*}
This is injective because in particular it maps $T_X$ to the $7$-dimensional irreducible subspace of $\Lambda^2$. It can be related to the bundle $W$ defined over the submanifold $\Sigma$ as follows. We consider the restriction
\begin{eqnarray*}
F:T_\Sigma&\to& \Lambda^2W^*\\
Y&\to& (Y\hook\varphi)|_W.
\end{eqnarray*}
That this is still injective is clear by taking $Y=e_1$ to be any unit vector in $T_\Sigma$. This is sufficient because of the transitive action of $G_2$ on the relevent space. Then according to the above description for $\varphi$, $F(Y)=\eta_1^-=\varepsilon^{45}-\varepsilon^{67}$ which is clearly non-zero on $W$. 

The spaces $T_\Sigma$ and $W$ are equipped with endomorphisms $\J$ with square minus the identity. We tensor these spaces with the complex numbers and extend $\J$ by complex linearity. We then obtain decompositions of them into $i$ and $-i$ eigenspaces. That is
\begin{eqnarray*}
T_\Comp=T\otimes \Comp =T^{1,0}\oplus T^{0,1},\ \ \ \ W_\Comp =W^{1,0}\oplus W^{0,1}.
\end{eqnarray*}
We extend the map $F$ by complex linearity.
\begin{prop}
The map $F$ sends $T^{1,0}$ isomorphically onto $\Lambda^{2,0}W^*$. Moreover, this map preserves the connections on the respective spaces. 
\end{prop}

 We will show that the image of $T^{1,0}$ is $\Lambda^{2,0}W^*$. This is again clear. Specifically, we can take the vector $Y=e_1-i\J e_1=e_1-ie_2\in T^{1,0}$ and $\sigma^1\wedge \sigma^2=(\varepsilon^4+i\varepsilon^7) \wedge(\varepsilon^6+i\varepsilon^5)\in\Lambda^{2,0}W^*$. Then,
\begin{eqnarray*}
F(Y)=\eta^-_1-i\eta^-_2=-i\sigma^1\wedge\sigma^2.
\end{eqnarray*}
That is, $F$ is a smooth isomorphism of these vector bundles. This means that the only numerical invariant of the bundle $W$ can be calculated solely from the intrinsic topological data of the submanifold. That is, $c_1(W)=-c_1(\Lambda^2W^*)=-c_1(T_\Sigma)=2g-2$.
The derivative of $F$ is also easily calculated.
\begin{eqnarray*}
(\nabla F)(Y,\xi)&=&\nabla(F(Y,\xi))-F(\nabla^T Y,\xi)-F(Y,\nabla^W\xi)\\
&=&(\nabla^X\varphi)(Y,\xi)+\varphi(\nabla^XY,\xi)+\varphi(Y,\nabla^X\xi)-\varphi(\nabla^TY,\xi)-\varphi(Y,\nabla^W\xi)\\
&=& \varphi(B(\cdot,Y),\xi)+\varphi(Y,D^{W^\perp}\xi)\\
&=&0.
\end{eqnarray*}
in the first case, the term vanishes because the second fundamental form $B$ takes values in $W$ and $\varphi$ vanishes on $W$, and in the second case because $\varphi(a,b,c)$ vanishes if exactly $2$ of $a,b$ and $c$ lie in the associative plane $T_\Sigma+\eta$ and the third is in $W$. Hence, $F$ is parallel. It is thus a holomorphic isomorphism of the line bundles $T_\Sigma$ and $\Lambda^2W^*$. 

\section{Examples} The first complete examples of manifolds with $G_2$ holonomy were constructed by Bryant and Salamon \cite{BS}. In the first case, they considered the total space of the bundle $\Lambda^-_M$ of anti-self-dual $2$-forms on a self-dual Einstein $4$-manifold $M$ of positive scalar curvature. This space can be expressed as the quotient $\mathcal{F}\times_\rho \text{Im}(\Ham)$ where $\mathcal{F}$ is the orthonormal coframe bundle on $M$. The torsion-free $G_2$ structure on $\Lambda^-_M$ is given as $\varphi=f\gamma_1+g\gamma_2$ where $f$ and $g$ are functions on $\Lambda^-_M$ and $\gamma_1$ and $\gamma_2$ are $3$-forms on $\mathcal{F}\times \text{Im}(\Ham)$ that are $SO(4)$ invariant and vanish on the kernel of the projection to $T(\Lambda^-_M)$. 

The zero section of $\Lambda^-_M$, considered as a submanifold, is totally geodesic for this metric. Along this subset, $f=g\equiv 1$ and $\varphi$ can be expressed as $\varphi =\varepsilon^{123}+\varepsilon^i\wedge\eta^-_i$. Here $\varepsilon^i=dx^i$ where $x^i$ give linear coordinates on the fibres of $\Lambda^-_M$, and $\{\eta^-_i\}$ form an orthogonal basis of anti-self-dual forms on $M$. 

Let $\Sigma\subseteq M$ be a  submanifold of $M$. Then, according to this expression for $\varphi$, the unit vector field $\eta$ to the submanifold is given as the vector dual to the covector 
\begin{eqnarray*}
\eta^\flat=\eta^-_i(\xi)\varepsilon^i = \xi\hook\varphi
\end{eqnarray*}
which is to say that $\eta$ is orthogonal to $M$ and hence orthogonal to all values of the second fundamental form. In particular, submanifolds of $S^4$ and $\CP^2$, considered as submanifolds of the total spaces of bundles,  satisfy the conditions of Theorem \ref{thm:holom-min}.

The other complete example given by Bryant and Salamon involved the total space of the spinor bundle on $S^3$. Similar considerations to those above show that submanifolds of the zero section of this bundle do not satisfy the required condition. The second fundamental form takes values colinear with $\eta$.

Another important example of a manifold with holonomy contained in $G_2$ is $X=M\times S^1$ where $(M,g,J,\Omega)$ is a Calabi-Yau manifold of real dimension $6$. The product $G_2$ structure is given by the $3$-form
\begin{eqnarray*}
\varphi =\text{Re}(\Omega) -dt\wedge\omega.
\end{eqnarray*}
$\Omega$ here is a parallel $(3,0)$-form on $M$ and $\omega$ is the K{\"a}hler form. For any complex curve $\Sigma\subseteq M$, considered as a submanifold of $M\times S^1$, the second fundamental form takes values tangent to $M$ but the distinguished unit normal vector field $\eta$ is given by 
\begin{eqnarray*}
\eta=(\xi\hook\varphi)^\#=\partial_t.
\end{eqnarray*}
Complex curves in $M\times\{t_0\}$ therefore satisfy the hypotheses of Theorem \ref{thm:holom-min} . However, arbitrary minimal submanifolds of $M$ do not necessarily.

\bibliographystyle{amsalpha}

\begin{thebibliography}{EFK}

\bibitem[B]{B} R.L. Bryant, {\em Conformal and minimal immersions of compact surfaces into the $4$-sphere}, J. Differential Geometry {\bf 17} (1982), 455-473.

\bibitem[BS]{BS} R.L. Bryant and S.M. Salamon, {\em On the construction of some complete metrics with exception holonomy}, Duke Math. J.  {\bf 58} (1989), 829-850.


\bibitem[ES]{ES} J. Eells and S. Salamon, {\em Twistorial constructions of harmonic maps into four-manifolds}, 
Ann. Scuola Norm. Sup. Pisa Cl. Sci. {\bf 12} (1985), 589-640.

\bibitem[G]{G} P. Gauduchon, {\em Pseudo-immersions superminimales d'une surface de Riemann dans une vari{\'e}t{\'e} riemannienne de dimension $4$}, Bull. Soc. Math. France {\bf 114} (1986), 447-508.

\bibitem[HL]{HL} R. Harvey and H.B. Lawson, Jr., {\em Calibrated geometries}, Acta Math.  {\bf 148} (1982), 47-157.

\bibitem[K]{K1} S. Karigiannis, {\em Some notes on $G_2$ and $Spin(7)$ geometry}, arXiv:math/0608618v3.

\bibitem[L]{L} F. Labourie, {\em Probl{\`e}me de Minkowski et surfaces {\`a} courbure constante dans les vari{\'e}t{\'e}s hyperboliques}, Bull. Soc. Math. France {\bf 119} (1991), 307-325.


\bibitem[S]{S} G. Smith, {\em Pointed $k$-surfaces}, Bull. Soc. Math. France {\bf 134} (2006), 509-557.

\bibitem[V]{V} M. Verbitsky, {\em A $CR$ twistor space of a $G_2$-manifold}, 	arXiv:1003.3170v2.

\end{thebibliography}

\end{document}